%From tadams@math.unc.eduTue Aug 13 11:00:44 1996
%Date: Fri, 9 Aug 1996 12:42:48 -0400 (EDT)
%From: Terry Adams <tadams@math.unc.edu>
%To: Karl Petersen <petersen@math.unc.edu>
%Cc: tadams@math.ohio-state.edu
%Subject: Re: Aug. 6 draft
\input amstex
\magnification=1200
\documentstyle{amsppt}
\parindent=1em
\nologo

\refstyle{C}
\loadeusm
\document

\def\colon{{:}\;}
\def\|{|\;}

\topmatter
\title Binomial-coefficient Multiples of Irrationals\endtitle
\author Terrence M. Adams and Karl E. Petersen\endauthor

\address Department of Mathematics, CB 3250 Phillips Hall, University of
North Carolina, Chapel Hill, NC 27599-3250, 
USA\endaddress
\curraddr
Department of Mathematics, Ohio State University, Columbus, OH 43210,
USA
\endcurraddr 
\email tadams\@math.ohio-state.edu \endemail

\address Department of Mathematics, CB 3250 Phillips Hall, University of
North Carolina, Chapel Hill, NC 27599-3250 
USA\endaddress 
\email petersen\@math.unc.edu \endemail 
\subjclass Primary 28D05, 28D99.\endsubjclass

%\date August 13, 1996\enddate
\abstract 
Denote by $x$ a random infinite path in the graph of Pascal's triangle
(left and right turns are selected independently with fixed probabilities)
and by $d_n(x)$ the binomial coefficient at the $n$'th level along the
path $x$. Then for a dense $G_{\delta}$ set of $\theta$ in the unit
interval, $\{d_n(x)\theta \}$ is almost surely dense but not uniformly
distributed modulo 1.
\endabstract 

%\thanks \endthanks

\keywords uniform distribution modulo 1, Pascal adic transformation,
measure-preserving transformation, eigenvalue, weak mixing\endkeywords 

\endtopmatter

\head  1. Introduction \endhead  

The {\it Pascal graph} is the directed infinite planar graph with vertices
$(n,k)$, for $n=0,1,\dots$ and $k=0,\dots ,n$ and two edges coming out
of each vertex $(n,k)$, one 
to $(n+1,k)$ and one to $(n+1,k+1)$. Labeling edges of the first kind by $0$
and of the second kind by $1$ produces a natural correspondence
between infinite sequences $x \in \Omega = \{ 0,1 \} ^{\Bbb N}$ and infinite
paths in the Pascal graph which start at the root vertex $(0,0)$. 
We will denote by $d_n(x)$ the binomial coefficient $C(n, k_n(x))$
found at the $n$'th vertex of $x$, if the Pascal graph is superimposed
on the Pascal triangle.
The {\it Pascal adic transformation} on the space $X$ of infinite
paths (see \cite{9, 11, 12, 6}) corresponds to the map $T: \Omega \to
\Omega$ given by $T(1^p0^q01\dots) = 0^q1^p10\dots$
for $p,q \geq 0$. Vershik \cite{9} noted that the invariant ergodic 
measures for this map are exactly the  
Bernoulli measures 
$\mu_{\alpha}=\Cal B(\alpha, 1 - \alpha)$ on $\Omega$ and 
conjectured \cite{10} that they are weakly mixing.
It was noted in \cite{6} that if,
for a fixed Bernoulli measure on $\Omega$, $\lambda$ is an eigenvalue
of $T$, then $\lambda ^{d_n(x)} \to 1$ for a.e. $x$.  The question of whether
or not there exist such $\lambda$, and its variants concerning the
distribution of the points $\lambda ^{d_n(x)}$ on the unit circle for
typical $x$ or indeed for all $x$, some of them also
mentioned in \cite{6}, lead to the study of the distribution modulo 1
of binomial-coefficient multiples of irrationals; answering many of these
questions seems to demand deeper understanding of the divisibility
properties of binomial coefficients than we have at present.

%\medskip
While we are not yet able to answer the question of weak mixing for
the Pascal adic transformation, we do have some progress on related
questions. First, we note that if $x$ is a path in the Pascal graph
which tracks a line of a fixed slope $\alpha$, then the proportion
of $j$'s, $0 \leq j \leq n$,  for which $d_j(x)$ is divisible by a 
fixed prime $q$
tends to $1$ as $n \to \infty$. 
Using this, we construct an uncountable set
of $\lambda$ on the unit circle such that for a.e. path $x$ in the
Pascal graph (with respect to a fixed Bernoulli measure) the points
$\lambda ^{d_n(x)}$ are not uniformly distributed on the circle, since
asymptotically too large a fraction of them are near $1$. Thus these
points might be candidates for eigenvalues of $T$---but we construct
such $\lambda$ for which the $\lambda ^{d_n(x)}$ are dense. We also list
several further observations, questions, and conjectures 
about the distribution of these points; perhaps the strongest conjecture
(also mentioned in \cite{6}) is the following:
If $\lambda \in \Bbb C$ and there exists a path $x$ in the
Pascal graph for which $\lambda ^{d_n(x)} \to 1$, then $\lambda = 1$.
%\medskip

The second author gratefully acknowledges the support of the Erwin
Schr\" odinger Institute, Vienna, where part of this research was conducted.
\bigskip

\head  2. Intersections of lines with Sierpinski's gasket\endhead

Let $E$ denote the triangle with vertices $(0,0)$, $(1,0)$, and 
$(0,1)$.  In this section we construct Sierpinski's gasket as a 
subset of $E$.   Then we prove that every straight 
line path (with slope strictly greater than 0) through Sierpinski's 
gasket intersects the gasket in a set with one-dimensional Lebesgue 
measure 0.  (Since Sierpinski's gasket has two-dimensional Lebesgue 
measure 0, Fubini's theorem guarantees that almost every line with a 
specified slope intersects Sierpinski's gasket with 
one-dimensional Lebesgue measure 0, but this is not sufficient for 
our purpose.)  We give two lemmas which set up the 
general proof.  Also we apply this to the generalized Sierpinski 
gasket defined at the end  of this section.

Given a natural number $n$, let $E_n^1$ denote the interior of the right 
triangle with vertices $((2^{n-1}-1)/2^{n-1},1/2^n)$, 
$((2^n-1)/2^n,0)$ and $((2^n-1)/2^n,1/2^n)$ and let 
$E_n^2$ be the interior of the triangle with vertices 
$(1/{2^n},(2^{n-1}-1)/{2^{n-1}})$, 
$(0,(2^n-1)/{2^n})$ and $(1/{2^n},(2^n-1)/{2^n})$.  We view 
the collection of $E_n^i$'s contained in the larger triangle $E$ with 
vertices $(0,0)$, $(1,0)$ and $(0,1)$.  The next two lemmas concern the 
intersection of the $E_n^i$'s with a 
straight line of slope $\gamma>0$ (where $\gamma$ corresponds to the 
invariant Bernoulli measure $\mu_{\alpha}$, with $\alpha = 1/(1+\gamma
)$ for the Pascal adic 
transformation).  

\proclaim{ Lemma 2.3}  Given $\alpha >0$ there exists $\eta=\eta 
(\alpha)>0$ such that: for any straight line $L$ of slope $\alpha$ which 
intersects the interior of $E$, there exists $n$ and $i$ such that
$$\frac{\mu (L\cap E_n^i)}{\mu (L\cap E)}\geq \eta.$$
\endproclaim 

Before we prove Lemma 2.3 we state and prove the preliminary Lemma 2.2.

\proclaim{ Lemma 2.2}  Given $\gamma >0$ there exist positive real 
numbers $\epsilon =\epsilon (\gamma)$ and $\eta =\eta (\gamma )$ such 
that for all $b\in [-\epsilon ,(1- \epsilon )/2]$ and 
$L_b=\{(x,\gamma (x-b)):x\in \Bbb R \}$ we have 
$$\frac{\mu (L_b\cap E_1)}{\mu (L_b\cap E)}\geq \eta.$$
\endproclaim

\demo{Proof}  First we find the intersection of $L_b$ with the hypotenuse 
of $E_1$.  Solving $y=\gamma (x-b)$ and $y=1/2 -x$ 
simultaneously, we obtain $x=(1/2 +\gamma b)/(1+\gamma )$ and 
$y={\gamma}(1/2 -b)/(1+\gamma )$.  Thus 
if we choose $\epsilon <1/(2\gamma )$ our point of intersection will 
have positive $x$ and $y$ coordinates.  Hence the function 
$f(b)={\mu (L_b\cap E_1 )}/{\mu (L_b\cap E)}$ 
is positive and continuous on the closed interval 
$[-\epsilon, (1- \epsilon )/2]$ and therefore achieves a 
positive minimum value $\eta$.\qed
\enddemo

\demo{Proof of Lemma 2.3}  Any line $L$ which intersects the interior of 
$E$ intersects either the line 
segment joining $(0,0)$ to $(1,0)$ or the line segment joining 
$(0,0)$ to $(0,1)$.  Without loss of generality let us consider 
line segments $L$ intersecting the segment joining $(0,0)$ to 
$(1,0)$. If we let $I=\{(x,0):0\leq x<1\}$, then $I\cap L\neq 
\emptyset$.  In this case we may focus on the sets $E_n=E_n^1$.

Choose $\epsilon >0$ as in Lemma 2.2.  (In particular $\epsilon 
<1/(2\gamma)$ will work.)  We may cover $I$ with subintervals 
$I_n=[(2^{n-1}-1)/{2^{n-1}}-{\epsilon}/{2^{n-1}}, 
(2^n-1)/{2^n}-{\epsilon}/{2^n})$;  so we have 
$I=\bigcup_{n=1}^{\infty} I_n$.  For each $b\in I_n$ define $L_b=\{ 
(x,\gamma (x-b)): x\in \Bbb R \}$ and $f_n:I_n\to [0,1]$ as 
$f_n(b)={\mu (L_b\cap E_n)}/{\mu (L_b \cap E)}$.  The self-similarity 
properties of the triangles $E_n$ imply that 
each $f_n:I_n\to [0,1]$ is continuous,  
and they all have the same image.  Therefore by Lemma 2.2 there exists a 
single real number $\eta >0$ such that 
$$f_n(b)=\frac{\mu (L_b \cap E_n)}{\mu (L_b \cap E)} \geq \eta$$
for all positive integers $n$ and all $b\in I_n$.\qed
\enddemo

Now we construct Sierpinski's gasket as a closed nowhere dense 
subset of $E$.  We call triangles in the plane {\it lower
triangles} if we can label the vertices $(a_1,b_1)$, $(a_2,b_2)$ and 
$(a_3,b_3)$ so that the right angle is at $(a_2,b_2)$ and 
$a_2=\min\{a_1,a_3\}$ and $b_2=\min\{b_1,b_3\}$.  {\it Upper triangles} 
have right angle at $(a_2,b_2)$ with 
$a_2=\max\{a_1,a_3\}$ and $b_2=\max\{b_1,b_3\}$.  Note that given a lower 
triangle $R$ there is a unique upper triangle (inscribed in 
$R$) whose vertices are the midpoints of the sides of $R$.  This upper 
triangle is denoted $\Cal{U}(R)$; let 
$\Cal{L}(R)=\{R_1,R_2,R_3\}$ be the collection of lower triangles 
remaining when we extract $\Cal{U}(R)$ from $R$.  Also given a collection 
$\Cal{C}$ of lower triangles let 
$\Cal{U}(\Cal{C})=\{\Cal{U}(R):R\in \Cal{C}\}$ and let 
$$\Cal{L}(\Cal{C})=\bigcup_{R\in \Cal{C}}\Cal{L}(R).$$ 
We take the triangles in $\Cal{L}(R)$ to be closed.

The following proposition uses Lemma 2.3 to prove that $\mu (L\cap G)=0$ 
for any line $L$ with slope $\gamma >0$.  First note that Lemmas 2.2 and 
2.3 can be extended to any lower triangle playing the role of
the initial triangle $E$.  Also note that 
$E_n^i\in \Cal{U}(\Cal{L}^{n-1}(E))$ for positive integers $n$ and 
$i=1,2$.  This implies that each $R\in \Cal{L}^n (E)$ is disjoint from 
the interior of $E_n^i$.  Let
$$G_n=\bigcup_{R\in \Cal{L}^n(E)} R.$$
Then $G=\bigcap_{n=1}^{\infty} G_n$ is {\it Sierpinski's gasket}.

\proclaim{Proposition 2.4} If $\gamma >0$ and
$L=\{(x,\gamma x):x\in \Bbb R \}$, then $\mu (G\cap 
L)=\lim_{n\to \infty} \mu (G_n \cap L)=0$.
\endproclaim

\demo{Proof}  Choose $\eta =\eta (\gamma) >0$ as in Lemma 2.3.  
We construct inductively a sequence $n_j$ of natural numbers such 
that for all positive integers $j$ we have
$$\frac{\mu (L\cap G_{n_j})}{\mu (L\cap E)}\leq (1-\eta )^j.$$
For the primary case Lemma 2.3 ensures that there exists 
$E_{n_1}^{i_1}\in \Cal{U}(\Cal{L}^{n_1-1}(E))$ such that 
$\mu (L\cap E_{n_1}^{i_1})/{\mu (L\cap E)}\geq \eta$.  (Actually 
$n_1=1$.)  Hence ${\mu (L\cap G_{n_1})}/{\mu (L\cap E)}\leq 1-\eta .$

For the general case suppose that $n_k$ satisfies ${\mu (L\cap G_{n_k})}/{\mu 
(L\cap E)}\leq (1-\eta )^k.$  
Now $\Cal{L}^{n_k}(E)=\{ R_1,\dots ,R_p\}$ is composed of a finite number of 
closed triangles.  Let $\Cal{M}=\{m:1\leq m\leq p,\ \ R_m\cap L \neq 
\emptyset\}$.  Thus $\mu (R_m\cap L)=0$ for $m\notin \Cal{M}$.  For each 
$m\in \Cal{M}$, by Lemma 2.3 there exists a positive integer $e(m)$ and 
$E(m)\in \Cal{U}(\Cal{L}^{e(m)-1}(R_m))$ such that 
${\mu (L\cap E(m))}/{\mu (L\cap R_m)}\geq \eta$.  Let $e=\max_{m\in 
\Cal{M}}\{e(m)\}$, and let $n_{k+1}=n_k +e$.  Therefore we have
$$\align 
\mu (L\cap G_{n_{k+1}})&\leq (1-\eta )\mu (L\cap G_{n_k}) \\ 
&\leq (1-\eta )^{k+1}.\qed
\endalign$$
\enddemo

Now we define the {\it generalized Sierpinski gasket} and give the analogous 
proposition which may be proved by the same method.  Let $E$ be the closed 
triangle described above.  Given a positive integer $q$ the lines 
$y={p}/{q}$, $x={p}/{q}$ and $y={p}/{q}-x$ partition $E$ into 
$q^2$ triangles with ${q(q-1)}/2$ upper triangles and 
${q(q+1)}/2$ lower triangles.  Let $\Cal{U}_q(E)$ be the 
collection of upper triangles and $\Cal{L}_q(E)$ be the collection 
lower triangles.  (Take the members of $\Cal{U}_q(E)$ to be open 
and the members of $\Cal{L}_q(E)$ to be closed.)  Define 
$\Cal{U}_q(\Cal{C})$ and $\Cal{L}_q(\Cal{C})$ analogously for any 
collection $\Cal{C}$ of lower triangles.  We obtain
$$G_q^n=\bigcup_{R\in \Cal{L}_q^n(E)}R\ ,\ \ \ G_q=\bigcap_{n=1}^{\infty} 
G_q^n.$$

\proclaim{Proposition 2.5}  If $\gamma >0$ and 
$L=\{(x,\gamma x):x\in \Bbb R \}$, then for all positive 
integers $q$ $$\mu (G_q \cap L)=\lim_{n\to \infty} \mu (G_q^n \cap L)=0.$$
\endproclaim

\demo{Proof}
We explain  how to choose sets $E_n^i(q)$ analogous to the
$E_n^i$ above.  
Then the analogues of the previous lemmas and proposition 
follow in the same manner as before.

Given a lower triangle $R$ let $\Cal{L}_q(R)$ be the collection of 
two lower triangles: the top left, lower triangle from 
$\Cal{L}_q(R)$ and the bottom right, lower triangle from 
$\Cal{L}_q(R)$.  Similarly we define $\Cal{L}_q^n(\Cal{C})$.  The sets 
$E_n^i(q)$ are chosen in $\Cal{U}_q(\Cal{L}_q^{n-1} (E))$.
\qed
\enddemo

\head  3. 
Binomial coefficients modulo a prime along a random path in Pascal's 
triangle \endhead

If Pascal's triangle is reduced modulo a prime $q$, a well-known
self-similar pattern (which can be produced by a cellular automaton)
results; this is a consequence of Kummer's Carry Theorem \cite{4}
and the resulting formula of Lucas \cite{5}. The parts of the
triangle that correspond to the upper triangles (which form $(G_q)^c$)
(the `voids') consist of regions in which the binomial coefficients
are divisible by $q$. The following theorem says that
 since a line of slope $\gamma$ spends most
of its time outside of each $G_q^n$, a $\mu_\alpha$-typical path $x$,
which eventually approaches a line of slope $\gamma=(1-\alpha)/\alpha$ in 
Pascal's
triangle, spends most of its time on vertices which carry binomial
coefficients divisible by $q$.

\proclaim {Theorem 3.1} 
If $q$ is prime and $0<\alpha <1$, then for $\mu _\alpha$-almost all
$x$ we have
$$
\lim_{n \to \infty} \frac {1}{n}\sum _{j=0}^{n-1}  \vert e^{2 \pi i d_j(x)/q} - 1
\vert = 0.
$$
\endproclaim

\demo{Proof} 
Let $\epsilon > 0$. Using Proposition 2.5, choose $N$ so large that if
$n \geq N$ then $\mu(G_q^n \cap L_\alpha )< \epsilon$. Notice that
since the complement of $G_q^N$ is a union of finitely many triangles,
if we move $L_\alpha$ just a small amount we cannot decrease the
Lebesgue measure of its intersection with $(G_q^N)^c$ by very
much. Thus we may choose $\delta > 0$ and then 
a large-enough natural number $M$ such 
that if the part within our unit triangle $E$ of the band of width
$\delta$ about the line $L_\alpha$ is cut into $M$ equally-spaced chunks
by lines parallel to the hypotenuse of $E$, 
and if one point is chosen
from each of those chunks, then the proportion of those points which
are in $G_q^n$ is still less than $2\epsilon$.

Let $S_k(x)$ denote the number of $1$'s along the path $x$ (regarded
as a sequence in $\Omega=\{ 0,1 \} ^{\Bbb Z}$).
By the Ergodic Theorem, for $\mu_\alpha$-almost every $x$ 
there is $K=K(x)$ such that 
$$
\vert S_k(x) - k \alpha \vert < k \delta \quad \text{for all } k \geq
K .
$$
Choose a time $M > K/\delta$, and 
consider Pascal's triangle down to
that level, including subtriangles of rank up to $N$. When this part
of Pascal's triangle is scaled down to lie over our unit triangle $E$,
we see the subtriangles that form $(G_q^N)^c$, and the scaled-down path
$x$ lies entirely inside the band of width $\delta$ about $L_\gamma$,
the line of slope $\gamma = (1-\alpha)/\alpha$ in the Sierpinski gasket.
We have arranged  that the proportion of vertices of the scaled-down
path
which are in $(G_q^N)^c$ is at least $1-2\epsilon$, and hence the
proportion of vertices of the path $x$ at which the binomial
coefficients $d_j(x)$ are divisible by $q$ is at least $1-2\epsilon$.
\enddemo

\head  4. Main result: a construction of special irrationals
\endhead

In this section we use Theorem 3.1 to construct an uncountable dense
set 
(in fact a $G_\delta$) 
$\Lambda \subset [0,1)$ such that for each $\theta \in \Lambda$ and for 
$\mu_{\alpha}$-almost every $x\in X$ we have that the sequence 
$\{d_j(x)\theta \}$ is not uniformly distributed modulo 1.  In particular we 
obtain a result similar to Theorem 3.1, but with the limit replaced by 
lim\,inf.

\proclaim{Theorem 4.1} There exist a dense
$G_\delta$ set 
$\Lambda \subset [0,1)$ 
and a set of full $\mu_\alpha$-measure $Y \subset X$
so that for each $\theta \in \Lambda$ and $x \in Y$ we have

$$\liminf_{n\to \infty} \frac1{n} \sum_{j=0}^{n-1} \vert e^{2\pi 
id_j(x)\theta}-1\vert =0.$$     

\endproclaim
 
\demo{Proof}  Let $\{ \epsilon_n\}$ be a sequence of positive real numbers 
satisfying $\sum_{n=1}^{\infty} \epsilon_n < \infty$, 
and let $\{ q_n\}$ be a sequence of primes increasing to $\infty$.
We will produce 
sequences $R_n$ of natural 
numbers and $\delta_n>0$ so that if 
$$Y_n=\{x\in X:\frac1{R_n} \sum_{j=0}^{R_n-1}\vert e^{2\pi 
id_j(x){p}/{q_n}}-1\vert <\frac1{n} \ \ \text{for}\ \ p=0,1,\dots ,q_n-1 
\},$$ 

$$\Lambda_n=\{ \theta \in [0,1):\text{ there exists}\ \ p=0,1,\dots 
,q_n-1\ \ \text{with}\ \ \vert \theta -\frac{p}{q_n}\vert 
< \delta_n\},$$

$$\Lambda=\bigcap_{k=1}^{\infty} \bigcup_{n=k}^{\infty} 
\Lambda_n\ \ \text{, and}\ \ \ \ Y=\bigcup_{k=1}^{\infty} 
\bigcap_{n=k}^{\infty} Y_n,$$
then 

$$ \liminf_{n\to \infty} \frac1{n} \sum_{j=0}^{n-1} \vert e^{2\pi 
id_j(x)\theta}-1\vert =0$$
for all $\theta \in \Lambda$  and all $x\in Y$.

For each $n$, by Theorem 3.1 
we may choose $R_n$ so that $\mu_{\alpha}(Y_n)>1-\epsilon_n$.  
Then, since 
$$
\{\frac{1}{R_n} \sum_{j=0}^{R_n-1} \vert e^{2\pi id_j(x)\theta}-1\vert: x
\in Y_n\}
$$ 
is a finite collection of
continuous functions of $\theta$, we may choose $\delta_n>0$ so that
$$\frac{1}{R_n} \sum_{j=0}^{R_n-1} \vert e^{2\pi 
id_j(x)\theta}-1\vert <\frac{2}{n}
$$
 for all $x \in Y_n$ and all 
$\theta \in \Lambda_n$.
Then $\Lambda$ is a dense $G_{\delta}$ with the usual 
topology, and $\mu_{\alpha}(Y)=1$ because 
$\sum_{n=1}^{\infty} \epsilon_n$ 
converges.  

To verify the outcome of the theorem, first choose $\theta \in \Lambda$ and 
$x\in Y$.  Then there exists a sequence $n_m\to\infty$ such that 
$\theta 
\in \Lambda_{n_m}$ for all positive integers $m$.  Also there exists $k$ 
such that $x\in Y_n$ for $n\geq k$.  Hence for $n_m \geq k$ we have 
$$\frac{1}{R_{n_m}} \sum_{j=0}^{R_{n_m}-1} \vert e^{2\pi 
id_j(x)\theta}-1\vert <\frac{2}{n_m}.\qed$$
\enddemo

\head 5. Density without uniform distribution \endhead 

In the previous section we constructed a dense $G_{\delta}$ set $\Lambda 
\subset 
[0,1)$ such that for each $\theta \in \Lambda$, $\{d_j(x)\theta\}$ is
not 
uniformly distributed 
modulo 1 for $\mu_{\alpha}$-almost every $x\in X$.  Those $\theta \in 
\Lambda$, which are irrational, remain candidates for eigenvalues of the 
Pascal adic transformation.  However in this section we will show that if 
the sequence $\{ \delta_n\}$ converges to zero 
sufficiently fast, then $\{ d_j(x)\theta \}$ is dense modulo 1 
for each $\theta 
\in \Lambda$ and for $\mu_{\alpha}$-almost every $x\in X$.  This excludes 
these $\theta$ as eigenvalues for the Pascal adic.  

We begin by considering Pascal's triangle modulo a prime $q$.  Recall 
that  for $n\in \Bbb N$ and $1\leq k\leq q^n-1$, we have 
$C(q^n,k)\cong_{q} 0$.  
Hence $C(q^n-1,k)\cong_{q}(-1)^k$ for $0\leq k\leq q^n-1$, which gives a 
`blocking line' on the triangle.  It is this `blocking line' that
yields total ergodicity of the Pascal adic.  
In Lemma 5.1 we note that among the binomial coefficients 
in the row numbered $q^n-2$ one can find all the congruence classes
modulo $q$, and in fact they appear in a regular way. This allows us
to show that along a random path in the triangle 
a hit of congruence class $r$ at level $q^m$ and of congruence class
$p$ at level $q^n$ are approximately independent if $m$ and $n$ are
far apart, and therefore with probability $1$  
no congruence class
modulo $q$ can be avoided. Consequently, if $\theta$ is very well
approximated by rationals $p_n/q_n$, then $\{ d_j(x)\theta \}$ must be
dense modulo $1$.

\proclaim{Lemma 5.1}  Let $q$ be prime and $n$ a natural number.  Then 
for $k=0,\dots ,q^n-2$ we have the following formula:
$$C(q^n-2,k)\cong_{q} (-1)^k(k+1). \tag 1$$
Moreover, for natural numbers $k$ and $p$ satisfying $0\leq p\leq q-1$ and 
$0\leq k\leq q^n-2q-1$ the set 
$$\{ i:k\leq i\leq k+2q-1,C(q^n-2,i)\cong_{q}p\} \tag 2$$
has exactly two elements. 
\endproclaim

\demo{Proof} First we derive the formula inductively.  The primary case 
is trivial: $C(q^n-2,0)=1=(-1)^0(1)$.  Assume the formula holds for 
$k=\ell -1$.  Thus for $k=\ell$ we have 
$$\align 
C(q^n-2,\ell )&= C(q^n-1,\ell )-C(q^n-2, \ell-1) \\ 
&\cong_{q} (-1)^{\ell} -(-1)^{\ell -1} \ell \\ 
&=(-1)^{\ell}[1+\ell ].
\endalign$$

Now for the second part of the lemma, we note that our formula gives the 
following: 
$$
C(q^n-2,i+2)\cong_{q} \cases C(q^n-2,i) + 2 &\text{if $i$ even}\\ 
C(q^n-2,i)-2 &\text{if $i$ odd}. \endcases 
$$
Hence if $q=2$ we obtain $C(q^n-2,k)\cong_{q} 1$ for $k$ even and 
$C(q^n-2,k)\cong_{q} 0$ for $k$ odd.  For $q\neq 2$ we have that each of the 
maps $j\mapsto C(q^n-2,k+2j)$ mod q and $j\mapsto C(q^n-2,k+1+2j)$ mod q gives 
a bijection of $\{ 0,\dots ,q-1\}$. \qed 
\enddemo 

\proclaim{Lemma 5.2}  Suppose that for each $n\in \Bbb N$, we have a 
unimodal distribution $f_n$ on the set $\{0,\dots ,q^n-2\}$.  
For each $p$ satisfying 
$0\leq 
p\leq q-1$, let 
$$M_p=\{m:0\leq m\leq q^n-2, C(q^n-2,m)\cong_{q} p\} .$$
If 
$\lim_{n\to \infty} \max\{f_n(m):0\leq m\leq q^n-2\} =0$, then 
$$\lim_{n\to \infty} \sum_{m\in M_p}f_n(m) = \frac1{q}.$$
\endproclaim

\demo{Proof}  We will show that 
$$\lim_{n\to \infty} (\sum_{m\in M_p}f_n(m) - \sum_{m\in M_r}f_n(m))=0$$
for all $p$ and $r$.  Without loss of generality, assume 
that 
$$\sum_{m\in M_p}f_n(m) \geq \sum_{m\in M_r}f_n(m).$$
Partition 
$\{0,\dots ,q^n-2\}$ into subintervals of $2q$ consecutive numbers with one 
remaining subinterval of at most $2q$ consecutive numbers.  Discard the 
subinterval which contains the peak of the distribution $f_n$, as well 
as its two adjacent subintervals, from the set $M_p$.  Call the 
remaining set $M_p^{\ast}$.  Now for each $m^{\ast} \in M_p^{\ast}$ there 
exists $m \in M_r$ in the next interval of length $2q$ towards the
peak of $f_n$ 
such that $f_n(m)\geq f_n(m^{\ast})$. 
Therefore
$$
\sum_{m\in M_r}f_n(m) \geq \sum_{m\in M_p^*}f_n(m) 
\geq \sum_{m\in M_p}f_n(m) -
6\max\{f_n(m):0\leq m\leq q^n-2\}.
$$
\qed 
\enddemo 

Let $F_n(q,p)$ be the set of paths which pass through a vertex $(q^n-2,k)$ 
satisfying $C(q^n-2,k)\cong_{q} p$.  Lemma 5.2 implies that the 
conditional probability of the set $F_n(q,p)$, given that the path passes 
through a fixed vertex, converges to $1/q$ as $n\to \infty$.  
Therefore for each $m\in \Bbb N$
$$\lim_{n\to \infty} \mu_{\alpha}(F_n(q,p)\cap F_m(q,r))=\frac1{q}\mu_\alpha(F_m(q,r)). \tag 3$$
A standard Hilbert space argument of A. R\'{e}nyi and P. 
R\'{e}v\'{e}sz \cite{8} implies that 
$F_n(q,p)$, $n=1,2,\dots$ is a mixing sequence of sets.  In particular 
we have Lemma 5.3, which says $F_n(q,p)$ `sweeps out'.  Finally we prove 
Theorem 5.4 using an approximation technique similar to that used in 
the previous section.  

\proclaim{Lemma 5.3}  For $0\leq p< q$ with $q$ prime and $0<\alpha <1$,
$$\mu_{\alpha}(\bigcup_{n=1}^{\infty} F_n(q,p))=1.$$
\endproclaim

\proclaim{Theorem 5.4}  There exist a dense $G_\delta$
set 
$\Lambda 
\subset [0,1)$ 
and a set of full $\mu_\alpha$-measure $Y \subset X$
so that for each $\theta \in \Lambda$ and $x \in Y$
the set 
$\{ e^{2\pi id_j(x)\theta}\colon j\in \Bbb 
N \}$ is dense (but not uniformly distributed) in $S^1$.
\endproclaim

\demo{Proof}  Let $\{\epsilon_n\}$ be a sequence of positive real numbers 
satisfying $\sum_{n=1}^{\infty} \epsilon_n < \infty$,
and let $\{ q_n\}$ be a sequence of primes increasing to $\infty$.
We will produce 
sequences $R_n$ of natural 
numbers and $\delta_n>0$ so that if 
$$Y_n=\bigcap_{p=0}^{q_n-1}\bigcup_{j=1}^{R_n}F_j(q_n,p),$$

$$\Lambda_n=\{ \theta \in [0,1):\text{ there exists}\ \ p_n=0,1,\dots 
,q_n-1\ \ \text{with}\ \ \vert \theta -\frac{p_n}{q_n}\vert 
< \delta_n\},$$

$$\Lambda=\bigcap_{k=1}^{\infty} \bigcup_{n=k}^{\infty} 
\Lambda_n\ \ \text{, and}\ \ \ \ Y=\bigcup_{k=1}^{\infty} 
\bigcap_{n=k}^{\infty} Y_n,$$
then  $\{ e^{2\pi id_j(x)\theta}\colon j\in \Bbb 
N \}$ is dense but not uniformly distributed 
for all $\theta \in \Lambda$ and all $x\in Y$.

For each $n$, 
by Lemma 5.3 we may 
choose $R_n$ so that $\mu_{\alpha}(Y_n)>1-\epsilon_n$.  Then choose 
$$\delta_n<\frac{1}{n\,C(q_n^{R_n}-2,(q_n^{R_n}-2)/2)}.$$
As before, if $\theta \in \Lambda$ and $x \in Y$, then we can find
arbitrarily large $n$ such that $\theta \in \Lambda_n$ and $x \in
Y_n$. Now $\theta$ is very well approximated by a rational $p_n/q_n$,
and as $p$ runs through the congruence classes modulo $q_n$, the
points $pp_n/q_n$ are $1/q_n$-dense modulo $1$. Further, for each
congruence class $p$ modulo $q_n$ there is $j=1,\dots ,R_n$ such that
at level $s=q_n^j$ the path $x$ has its binomial coefficient
$d_s(x)$ hit that congruence class. Since $\delta_n$ has been
chosen so small that all the points $d_s(x)\theta$ under consideration
are very close to the points $d_s(x)p_n/q_n$, and the latter are
$1/q_n$-dense, we are done.
\qed \enddemo

\head 6. Questions and Conjectures \endhead

1. {\it Conjecture} \cite{10, 6}: 
For each Bernoulli measure $\mu _\alpha$, the Pascal adic
transformation $T$ is weakly mixing.  This would follow if 
one could show that $\lambda^{d_n(x)} 
\to 1$ for a.e. $x$ with respect to $\mu _\alpha$ implies $\lambda =1$.

%\medskip
2. {\it Conjecture}: 
If there is a path $x$ such that $\lambda^{d_n(x)} \to 1$, then
$\lambda = 1$.

%\medskip
3. Does there exist any $\lambda$ in the unit circle such that
$\{\lambda^{d_n(x)}\}$ is uniformly distributed in the circle
{\it for every $x$} (except for the two paths down the edges)? For
such a $\lambda$, the skew-product transformation
$$
S(z_1,z_2,z_3,\dots) = (\lambda z_1, z_1z_2, z_2z_3, \dots)
$$
on the infinite torus, known to be uniquely ergodic by results of
Weyl \cite{13, 14}, Furstenberg \cite{1}, Hahn \cite{3}, and 
Postnikov \cite{7}, would have the very strong
property that we would see a uniformly distributed sequence 
$\{(S^jz)_k\}$ not only when we
looked in a fixed coordinate $k$ at the orbit of a point $z$, but also when we
allowed our view to shift one place to the right from time to time: 
$\{(S^jz)_{k_j}\}$ would be uniformly distributed for each $z$ and each
choice of $\{k_j\} \subset \Bbb N$ with $k_{j+1}-k_j \in \{0,1\}$ for
each $j$. (The $\lambda$ that we construct above are at another
extreme from this property.)

%\medskip
4. From another theorem of Weyl and Tonelli's Theorem it follows that
for almost every $\lambda$ the sequence $\{\lambda^{d_n(x)}\}$ is
uniformly distributed for a.e. $x$, with respect to each Bernoulli
measure $\mu_\alpha$. For which $\lambda$ does this hold? Similarly,
what paths $x$ have the property that this sequence is uniformly
distributed for each $\lambda$ that is not a root of unity? (By Weyl's
Theorem, this is the case for each path $x$ that is eventually
diagonal.)

%\medskip
5. Studies like these on divisibility of binomial coefficients by
primes suggest questions on simultaneous divisibility by several
primes. For example, thinking about the central path in Pascal's
triangle and divisibility by 2 and 3 leads to the
following {\it Conjecture}:
The only solutions in nonnegative integers $r$ and distinct $s_1, \dots
,s_m$ of an equation
$$
2^r=3^{s_1}+\dots +3^{s_m}
$$
are $1=1$, $4=1+3$, and $256=1+3+9+243$.
We thank Charles Giffen for pointing out that this conjecture was
already made by Erd\" os---see \cite{2}.
More generally, if for a prime $p$ we define an integer Cantor set
$H(p)$ to consist of all those expansions base $p$ with coefficients
in the interval $[0,p/2]$ (or subject to some other restriction),  
is $H(p_1) \cap \dots \cap H(p_n)$ typically finite?

\enddocument